# Multiresolution Solution of Burgers Equation with B-spline Wavelet Basis

**S. Hadi Seyedi**[a]

[a]Mechanical Engineering Department, Wayne State University, Detroit, Michigan 48202, USA

**Abstract:**

This paper represents a mixed numerical method for the multi-resolution solution of non-linear partial differential equations based on B-Spline wavelets. The method is based on a second-order finite difference formula combined with the collocation method which uses the wavelet basis and applied to the Burgers equation. Performance and accuracy of the numerical solutions are studied using three standard test cases with Dirichlet and Neumann boundary conditions. Comparing the results with other methods such as the fully implicit finite difference method and mixed finite difference/boundary element method shows a greater accuracy while using smaller number of basis functions.

## 1 Introduction

Burgers equation is one of the most well-known examples of nonlinear partial differential equations. This equation is important in a variety of applications such as gas dynamics and traffic flow [1], sound waves in a viscous medium and waves in fluid filled viscous elastic tubes. However, the most notable application can be found in simplification of the Navier-Stokes equation which simulates fluid dynamics problems such as flow through a shock wave traveling in a viscous fluid.

To demonstrate some specifications of turbulent fluid in a channel that is caused by effects of convection and diffusion terms, Burgers offered an equation [2, 3]. The Burgers equation established a coupling between diffusion term $u_{xx}$ with diffusive effects, and the convective term $uu_x$ with propagation effects. Therefore, solutions for this equation demonstrate a good balance between advection and diffusion terms.

Since for many combination of initial and boundary conditions, analytical expressions are available [4, 5], this equation is an appropriate choice for testing and comparing numerous numerical techniques.

$$u_t + uu_x = \frac{1}{\text{Re}} u_{xx} \qquad (x,t) \in [a,b] \times [0,T] \qquad (1)$$

Many numerical methods are used in the literature for solving Burgers equation [6-8]. For instance, one can find Galerkin finite element method [9], Homotopy analysis method [10], Adomian–Pade technique [11], spectral collocation method [12, 13], Fourth-order finite difference method [14] and Galerkin method with interpolating scaling functions [15] in the literature.

Numerical methods used for solving differential equations usually fall into two distinct main classes, i.e. local and global. The finite-difference and finite-element methods are based on local arguments, whereas the spectral method is global in character. Finite difference and finite element are suitable for modeling complex geometries, whereas spectral methods can provide superior accuracy [16]. To combine the advantage of both finite difference/finite element and spectral methods, wavelet basis are a new choice. Moreover, results can be obtained in multi scale form.

Wavelet theory is a relatively new area in numerical research for solving differential equations. Wavelets can be separated into two distinguishable types, orthogonal and semi-orthogonal [17]. Most papers in this field use orthogonal wavelet basis which almost have infinite support, or non-symmetric properties. In contrast, semi-orthogonal wavelets have finite support, simple analytical expressions and symmetry properties which provide ideal properties for characterization of a function [18]. The method of weighted residuals, especially Galerkin and collocation techniques are common for solving differential equations using wavelets [19-24].

The method, which is used in this paper, is a mixed finite difference collocation method (MFDCM) which creates a semi-discrete form of Burgers equation by applying the collocation method with B-spline wavelet basis for the spatial part. A second order finite difference method is used for the temporal part to fully discretize the equation.

## 2 B-spline wavelets

B-spline wavelets are constructed from B-splines functions of order $m$. An $m$ th order B-spline function divides the support $[0,1]$ to $m$ segments. The dual function of this spline will have $2m-1$ segment in this interval. For any resolution $s$, the domain is divided to $2^s$ segments with the step size of $1/2^s$. Therefore, the lowest octave level (to have at least one inner wavelet) can be determined by [18].

$$2^s \geq 2m-1 \tag{2}$$

If spline cardinal functions are used directly as the basis for expanding a function, then some parts of scaling functions and their wavelets will be placed outside of the domain of interest [22]. A remedy to this, is to use a compactly supported basis function for finite domains. One will need boundary and scaling functions for this purpose. Scaling functions and wavelets suggested in [18, 22] are used here. Only $m = 2$ is used in formula 2 and leads to having linear spline wavelets.

The linear B-spline scaling functions are given by:

$$\phi_{j,k}(x) = \begin{cases} x_j - k & k \leq x_j < k+1 \\ 2-(x_j - k) & k+1 \leq x_j < k+2 \quad , k = 0,\ldots,2^j - 2 \\ 0 & o.w. \end{cases} \qquad (3)$$

where $j$ is the resolution level and $k$ is the transfer parameter. For the left and right-side boundary scaling functions, $\phi_{j,k}(x)$ is respectively calculated as follows:

$$\phi_{j,k}(x) = \begin{cases} 2-(x_j - k) & 0 \leq x_j < 1 \\ 0 & o.w. \end{cases} \quad , k = -1 \qquad (4)$$

$$\phi_{j,k}(x) = \begin{cases} x_j - k & k \leq x_j < k+1 \\ 0 & o.w. \end{cases} \quad , k = 2^j - 1$$

In Eqs. 3 and 4 the actual position parameter $x$ is related to $x_j$ by $x_j = 2^j x$.

B-spline wavelets are given by:

$$\psi_{j,k}(x) = \frac{1}{6} \begin{cases} x_j - k & k \leq x_j < k+\frac{1}{2} \\ 4 - 7(x_j - k) & k+\frac{1}{2} \leq x_j < k+1 \\ -19 + 16(x_j - k) & k+1 \leq x_j < k+\frac{3}{2} \\ 29 - 16(x_j - k) & k+\frac{3}{2} \leq x_j < k+2 \quad , k=0,\ldots,2^j - 3 \\ -17 + 7(x_j - k) & k+2 \leq x_j < k+\frac{5}{2} \\ 3 - (x_j - k) & k+\frac{5}{2} \leq x_j < k+3 \\ 0 & o.w. \end{cases} \qquad (5)$$

For the left and right-side boundary wavelets, $\psi_{j,k}(x)$ is respectively calculated using:

$$\psi_{j,k}(x) = \frac{1}{6}\begin{cases} -6+23x_j & 0 \le x_j < \frac{1}{2} \\ 14-17x_j & \frac{1}{2} \le x_j < 1 \\ -10+7x_j & 1 \le x_j < \frac{3}{2} \\ 2-x_j & \frac{3}{2} \le x_j < 2 \\ 0 & o.w. \end{cases}, k = -1 \quad (6)$$

$$\psi_{j,2^j-2}(x) = \frac{1}{6}\begin{cases} 2-(k+2-x_j) & k \le x_j < k+\frac{1}{2} \\ -10+7(k+2-x_j) & k+\frac{1}{2} \le x_j < k+1 \\ 14-17(k+2-x_j) & k+1 \le x_j < k+\frac{3}{2} \\ -6+23(k+2-x_j) & k+\frac{3}{2} \le x_j < k+2 \\ 0 & o.w. \end{cases}, k = 2^j - 2 \quad (7)$$

Figure 1 shows the linear spline scaling functions and first wavelets on $[0,1]$

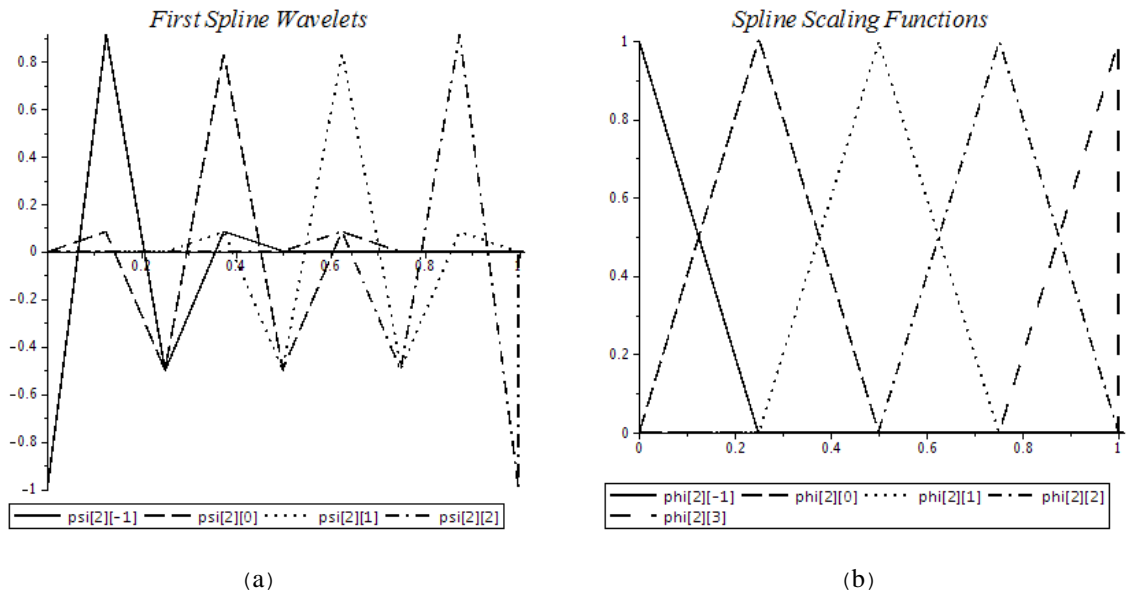

(a)                 (b)

**Figure 1.** Linear spline wavelets (a) and scaling functions (b) on $[0, 1]$

## 2.1 Approximation of a function on an interval

A function $f(x)$ can be approximated by $f_n(x)$ with B-spline wavelets over $[0,1]$

$$f(x) \approx f_n(x) = \sum_{k=-1}^{3} a_k \phi_{2,k} + \sum_{i=2}^{M} \sum_{j=-1}^{2^i-2} b_{i,j} \psi_{i,j} = C^T \Psi \tag{8}$$

Where $\phi_{2,k}$, $\psi_{i,j}$ are scaling and wavelet functions respectively. $C$ and $\Psi$ are vectors given by

$$C = \left[ c_{-1}, c_0, \ldots, c_3, d_{2,-1}, \ldots, d_{2,2}, d_{3,-1}, \ldots, d_{3,6}, \ldots, d_{M,-1}, d_{M,2^M-2} \right]^T \tag{9}$$

$$\Psi = \left[ \phi_{2,-1}, \phi_{2,0}, \ldots, \phi_{2,3}, \psi_{2,-1}, \ldots, \psi_{2,2}, \psi_{3,-1}, \ldots, \psi_{3,6}, \ldots, \psi_{M,-1}, \psi_{M,2^M-2} \right]^T$$

Since B-spline wavelets are semi-orthogonal, for finding $C_i$ in Eq. 8, one should use dual functions of related basis.

$$C_i = \int_0^1 f(x) \bar{\Psi}_M(x) dx \tag{10}$$

In which $\bar{\Psi}_M$ are dual functions of $\Psi_M$. These dual functions can be evaluated by linear combinations of $\Psi_M$ considering:

$$\int_0^1 \bar{\Psi}(x) \Psi^T(x) dx = I \tag{11}$$

in which $I$ is the identity matrix. Defining $P_M$ as

$$P_M = \int_0^1 \Psi_M(x) \Psi_M^T(x) dx \tag{12}$$

From Eqs. 11 and 12 and some manipulation, one can conclude:

$$\bar{\Psi}_M(x) = (P_M)^{-1} \Psi_M(x) \tag{13}$$

Thus, unknown coefficients $C_i$ can be evaluated.

## 2.2 The operational matrix of derivatives

Derivatives of the base vector can be expressed as

$$\Psi'(x) = \begin{bmatrix} \Psi'_1 \\ \Psi'_2 \\ \vdots \\ \Psi'_n \end{bmatrix} = [OD] \begin{bmatrix} \Psi_1 \\ \Psi_2 \\ \vdots \\ \Psi_n \end{bmatrix} \qquad (14)$$

Where $[OD]$ is a $2^M + 1 \times 2^M + 1$ matrix. For constructing matrix of derivatives $[OD]$, one can use orthogonality and semi-orthogonality properties. For example, for the first row of the matrix one can do as follows:

$$OD_{11} = \int_0^1 \Psi'_1(x)\bar{\Psi}_1(x)dx$$

$$OD_{12} = \int_0^1 \Psi'_1(x)\bar{\Psi}_2(x)dx \qquad (15)$$

$$\vdots$$

$$OD_{1n} = \int_0^1 \Psi'_1(x)\bar{\Psi}_n(x)dx$$

## 2.3 The mixed finite difference and collocation method (MFDCM) for Burgers equation

In this section, nonlinear Burgers equation will be solved on a bounded domain. To this aim, finite difference method is used for discretizing the PDE to a system of ordinary differential equations. Then these ODE equations will be solved by wavelet collocation method. Time discretization step is set to $\Delta t = 10^{-3}$ for all cases.

By discretizing Eq. 1 according to the following θ-weighted scheme one will have:

$$\frac{u(x, t + \Delta t) - u(x, t)}{\Delta t} + u(x,t)u_x(x,t) + \theta(-\frac{1}{\text{Re}} u_{xx}(x, t + \Delta t)) \qquad (16)$$
$$+ (1 - \theta)(-\frac{1}{\text{Re}} u_{xx}(x,t)) = 0$$

Where $0 \leq \theta \leq 1$ and $\Delta t$ is the time step size. By rearranging Eq. 16 and using Crank-Nicolson method with $\theta = 1/2$ for diffusion part, one can obtain:

$$u(x, t + \Delta t) - \frac{\Delta t}{2\text{Re}} u_{xx}(x, t + \Delta t) = u(x,t) + \frac{\Delta t}{2\text{Re}} u_{xx}(x,t) - \Delta t \cdot u(x,t)u_x(x,t) \qquad (17)$$

Using Eq. 8, the function $u(x,t)$ can be approximated as:

$$u(x,t) = C^T \Psi \tag{18}$$

Also, by using Eq. 14 we can write:

$$u_x(x,t) = u(x,t)[OD]\Psi \tag{19}$$
$$u_{xx}(x,t) = u(x,t)[OD]^2 \Psi$$

Replacing Eqs. 18 and 19 in Eq. 17 and using the notation $u^n = u(x,t^n)$ where $t^n = t^{n-1} + t.\Delta t$, the equation will be changed to:

$$C_{n+1}^T \Psi - \frac{\Delta t}{2\operatorname{Re}}(C_{n+1}^T OD^2 \Psi) = C_n^T \Psi + \frac{\Delta t}{2\operatorname{Re}}(C_n^T OD^2 \Psi) - \Delta t(C_n^T \Psi)(C_n^T OD\Psi) \tag{20}$$

By collocating points in Eq. 20 on $2^M - 1$ points at $x_i = \frac{j-1}{N_p - 1}$, $j = 2,...,N_p - 1$ ($N_p$ is the total number of the coefficients) a system of $2^M - 1$ equations will be obtained.

For the boundary conditions one can collocate the boundary points $a$, $b$ in the solution:

$$C_{n+1}^T \Psi(a) = B.C_1 \tag{21}$$
$$C_{n+1}^T \Psi(b) = B.C_2$$

Where $B.C_1$ and $B.C_2$ are the left and right boundary conditions, respectively. Now a system of $2^M + 1$ equations and $2^M + 1$ unknowns is available. Equations are in the form of $AC_{n+1} = b$, $n = 0,1,...$ and $C_1$ can easily be determined using initial condition as follows ($G_0(x)$ is the initial condition of the problem)

$$G_0(x) = C_0^T \Psi(x) \tag{22}$$

This system of equations can be solved to find $C_{n+1}$ in any time step $(n = 0,1,...)$

### 3  Numerical results

In this section, three test problems are studied to validate and check the performance of the proposed scheme with different types of boundary conditions.

### 3.1 Test case 1.

In this section the Burgers equation (Eq. 1) with Reynolds number equal to 1.0 and 10.0 are presented. The initial and boundary conditions are:

$$\frac{\partial u}{\partial t} + u \frac{\partial u}{\partial x} = \frac{1}{Re} \frac{\partial^2 u}{\partial x^2} \quad x \in [0,1], t \geq 0 \quad (23)$$

I.C. $u(x,0) = \sin(\pi x)$

B.C. $u(0,t) = u(1,t) = 0$

The exact solution for this test case is [5]:

$$u(x,t) = \frac{\frac{2\pi}{Re}\sum_{n=1}^{\infty} a_n \exp(-\frac{1}{Re}n^2\pi^2 t) n \sin(n\pi x)}{a_0 + \sum_{n=1}^{\infty} a_n \exp(-\frac{1}{Re}n^2\pi^2 t)\cos(n\pi x)} \quad (24)$$

Where the Fourier coefficients are calculated by:

$$a_0 = \int_0^1 \exp\{-(\frac{2\pi}{Re})^{-1}[1-\cos(\pi x)]\}dx \quad (25)$$

$$a_n = 2\int_0^1 \exp\{-(\frac{2\pi}{Re})^{-1}[1-\cos(\pi x)]\}\cos(n\pi x)dx, \quad n = 1,2,3,...$$

For testing the accuracy of the discussed method, a comparison with the exact solution and two numerical schemes, the fully implicit finite difference method (IFDM) [6] and the mixed finite difference and boundary element method (BEM) [7] is performed

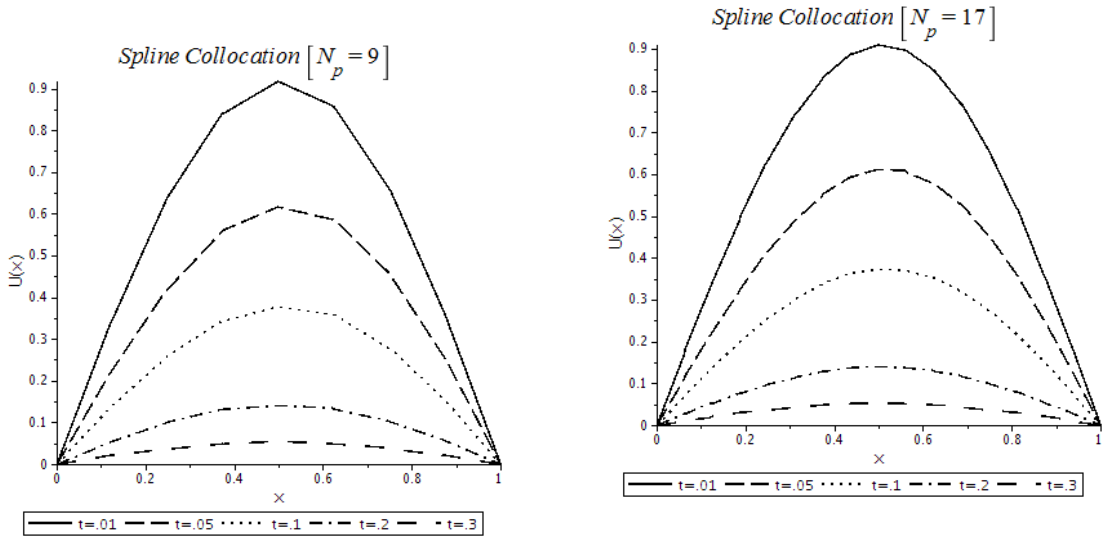

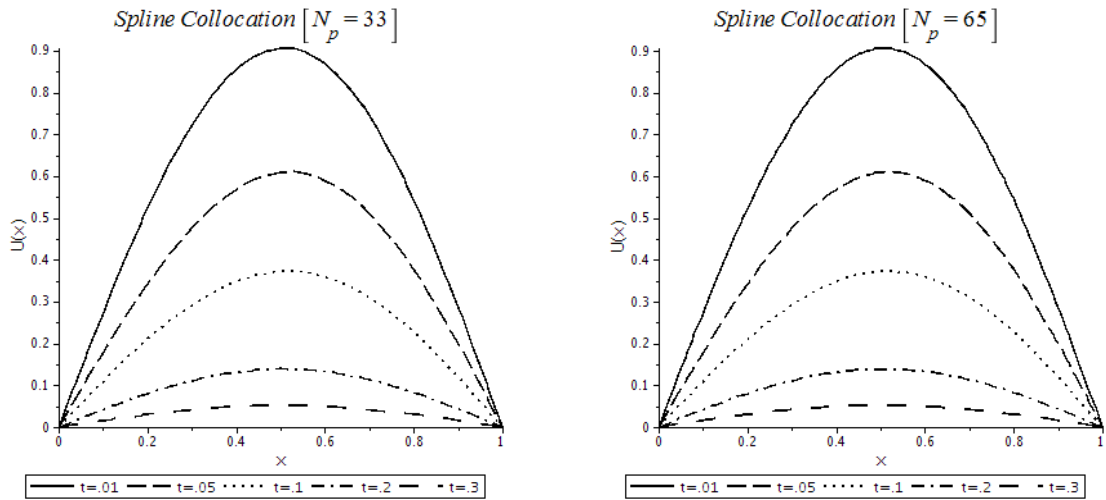

**Figure 2.** Numerical solutions in different scales for test case 1 with $Re = 1.0$

**Table 1.** Comparison of the results for test case 1 with $Re = 1.0$

| Method | time | x=.1 | x=.3 | x=.5 | x=.7 | x=.9 |
|---|---|---|---|---|---|---|
| IFDM[6] | 0.05 | 0.17832 | 0.47658 | 0.60984 | 0.51165 | 0.20006 |
| BEM[7] | | 0.17759 | 0.47531 | 0.60851 | 0.51050 | 0.19933 |
| EXACT | | 0.17803 | 0.47586 | 0.60907 | 0.51113 | 0.19989 |
| Current work ($N_p$=33) | | 0.17798 | 0.47558 | 0.60954 | 0.51109 | 0.20005 |
| Current work ($N_p$=65) | | 0.17792 | 0.47570 | 0.60918 | 0.51130 | 0.19999 |
| IFDM[6] | 0.1 | 0.11009 | 0.29335 | 0.37342 | 0.31144 | 0.12128 |
| BEM[7] | | 0.10931 | 0.29124 | 0.37070 | 0.30911 | 0.12031 |
| EXACT | | 0.10954 | 0.29190 | 0.37158 | 0.30991 | 0.12069 |
| Current work ($N_p$=33) | | 0.10952 | 0.29174 | 0.37187 | 0.30988 | 0.12078 |
| Current work ($N_p$=65) | | 0.10948 | 0.29181 | 0.37163 | 0.30999 | 0.12073 |
| IFDM[6] | 0.2 | 0.04273 | 0.11276 | 0.14120 | 0.11574 | 0.04457 |
| BEM[7] | | 0.04220 | 0.11044 | 0.13809 | 0.11322 | 0.04391 |
| EXACT | | 0.04193 | 0.11062 | 0.13847 | 0.11347 | 0.04369 |
| Current work ($N_p$=33) | | 0.04194 | 0.11059 | 0.13860 | 0.11345 | 0.04371 |
| Current work ($N_p$=65) | | 0.04192 | 0.11061 | 0.13849 | 0.11348 | 0.04369 |

Figure 2 shows the numerical results obtained in the current work for test case 1 with $Re = 1$, using 9, 17, 33 and 65 collocation points ($N_p$). It can be seen that the trend of the solution is identical, and the values are very close in all four cases even when using a low number of collocation points. To compare the results quantitatively, the value of the numerical solution obtained in the current work is given in Table 1, beside the exact solution and two other numerical solutions from the literature [6, 7]. It is evident from this table that by using 33 collocation points, the solution error is less than the other two methods. The IFDM solution was obtained using 100 grid points [6].

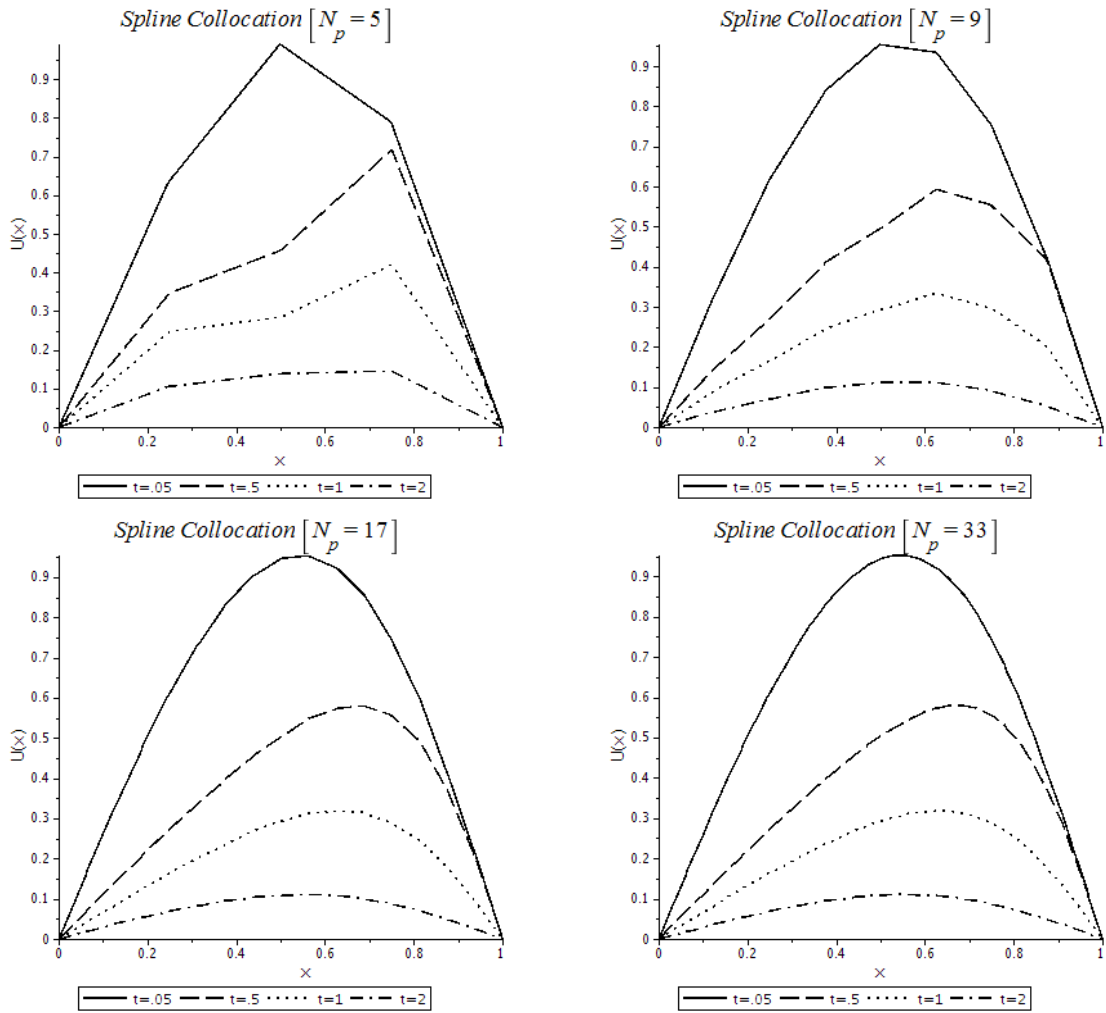

**Figure 3.** Numerical solutions in different scales for test case 1 with $Re = 10.0$

**Table 2.** Comparison of results for test case 1 with $Re = 10.0$

| Method | time | x=.1 | x=.3 | x=.5 | x=.7 | x=.9 |
|---|---|---|---|---|---|---|
| IFDM[6] | 0.5 | 0.11048 | 0.32367 | 0.50447 | 0.57664 | 0.30912 |
| BEM[7] |  | 0.10986 | 0.32191 | 0.50240 | 0.57514 | 0.30779 |
| EXACT |  | 0.10992 | 0.32219 | 0.50279 | 0.57585 | 0.30935 |
| Current work ($N_p$=33) |  | 0.10991 | 0.32207 | 0.50281 | 0.57542 | 0.30933 |
| IFDM[6] | 1 | 0.06689 | 0.19445 | 0.29448 | 0.31107 | 0.14769 |
| BEM[7] |  | 0.06644 | 0.19263 | 0.29139 | 0.30711 | 0.14507 |
| EXACT |  | 0.06632 | 0.19279 | 0.29192 | 0.30809 | 0.14607 |
| Current work ($N_p$=33) |  | 0.06631 | 0.19270 | 0.29194 | 0.30776 | 0.14605 |
| IFDM[6] | 2 | 0.02909 | 0.08044 | 0.10939 | 0.09838 | 0.04037 |
| BEM[7] |  | 0.02913 | 0.07951 | 0.10770 | 0.09663 | 0.03976 |
| EXACT |  | 0.02876 | 0.07946 | 0.10789 | 0.09685 | 0.03969 |
| Current work ($N_p$=33) |  | 0.02875 | 0.07941 | 0.10792 | 0.09676 | 0.03968 |

Figure 3 is the plot of the numerical solution for test case 1 with the $Re$ set to 10. As this figure shows, the convection effects become apparent in the solution and for the three cases of $N_p = 9, 17, 33$ the trend and the values of the obtained solutions are very similar. For quantitative comparison, the obtained numerical solution at certain locations are given in Table 2 vs the exact solution and the IFDM [6] and the BEM [7] solutions. The accuracy of the solution shows the same trend as with $Re = 1$

### 3.2 Test case 2.

For this test case the burgers equation is solved with the following initial and boundary conditions:

I.C.  $u(x,0) = 4x(1-x)$  (26)

The analytical solution for this case is the same as the previous one, but the Fourier coefficients are:

$$a_0 = \int_0^1 \exp[-x^2(\frac{3}{Re})^{-1}(3-2x)]dx \qquad (27)$$

$$a_n = 2\int_0^1 \exp[-x^2(\frac{3}{Re})^{-1}(3-2x)]\cos(n\pi x)dx, \quad n = 1,2,3,...$$

Again, in this test case, Reynolds number is set equal to 1.0 and 10.0. The results for this test problem are presented in Tables 3 and 4 and Figs 6 and 7.

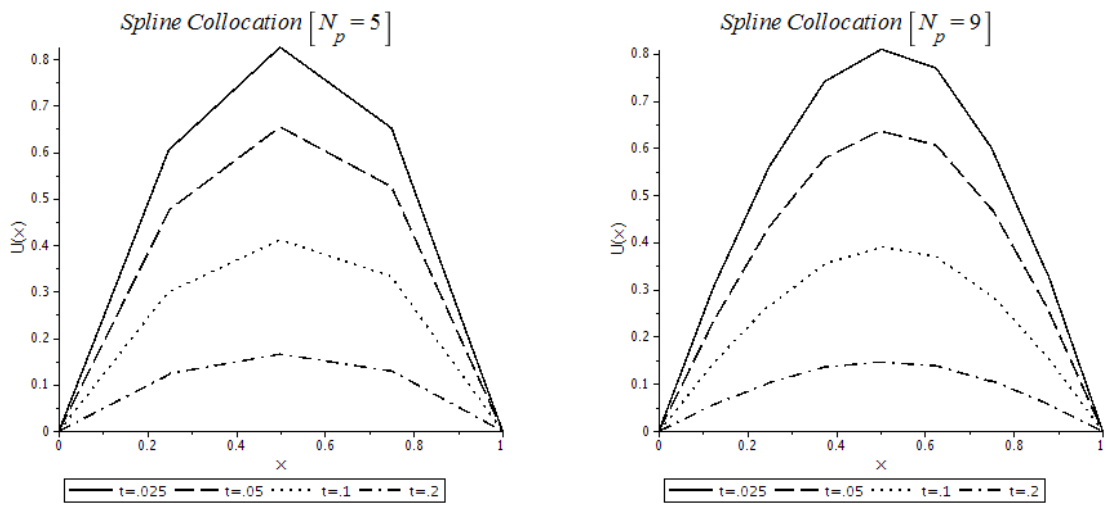

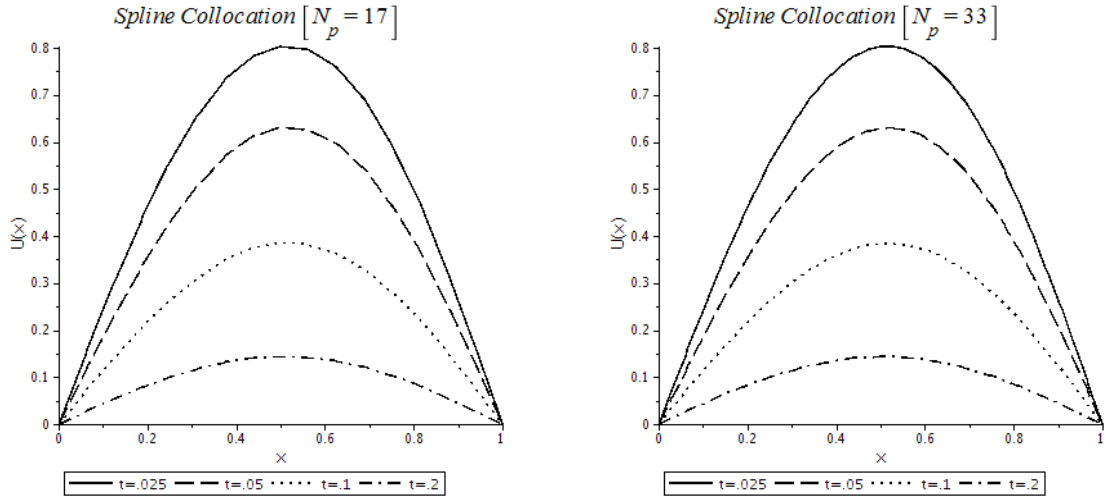

**Figure 4.** Numerical solution in different scales for test case 2 with $Re = 1.0$

Solutions obtained for this test case are similar to the test case 1. From a qualitative point of view, even with the lowest number of collocation points ($N_p = 5$), the evolution of the solution is similar to the one with the highest number of collocation points ($N_p = 33$), as can be seen in Figs. 4 and 5. Quantitative values are given in Tables 3 and 4 and the error of the numerical solutions relative to the exact solutions can be calculated. Table 5 gives the average of the relative error for the IFDM, the BEM and the current work with $Re = 1$ and $Re = 10$.

**Table 3.** Comparison of the results for test case 2 with $Re = 1.0$

| Method | time | x=.1 | x=.3 | x=.5 | x=.7 | x=.9 |
|---|---|---|---|---|---|---|
| IFDM[6] | 0.05 | 0.18423 | 0.49169 | 0.62884 | 0.52847 | 0.20712 |
| BEM[7] |  | 0.18347 | 0.49036 | 0.62749 | 0.52726 | 0.20632 |
| EXACT |  | 0.18389 | 0.49093 | 0.62808 | 0.52793 | 0.20690 |
| Current work ($N_p$=33) |  | 0.18385 | 0.49067 | 0.62858 | 0.52790 | 0.20707 |
| IFDM[6] | 0.1 | 0.11346 | 0.30248 | 0.38533 | 0.32165 | 0.12533 |
| BEM[7] |  | 0.11266 | 0.30031 | 0.38251 | 0.31925 | 0.12432 |
| EXACT |  | 0.11289 | 0.30097 | 0.38342 | 0.32007 | 0.12472 |
| Current work ($N_p$=33) |  | 0.11288 | 0.30082 | 0.38374 | 0.32005 | 0.12481 |
| IFDM[6] | 0.2 | 0.04407 | 0.11631 | 0.14570 | 0.11948 | 0.04602 |
| BEM[7] |  | 0.04353 | 0.11393 | 0.14250 | 0.11691 | 0.04535 |
| EXACT |  | 0.04324 | 0.11410 | 0.14289 | 0.11713 | 0.04511 |
| Current work ($N_p$=33) |  | 0.04325 | 0.11407 | 0.14302 | 0.11712 | 0.04514 |

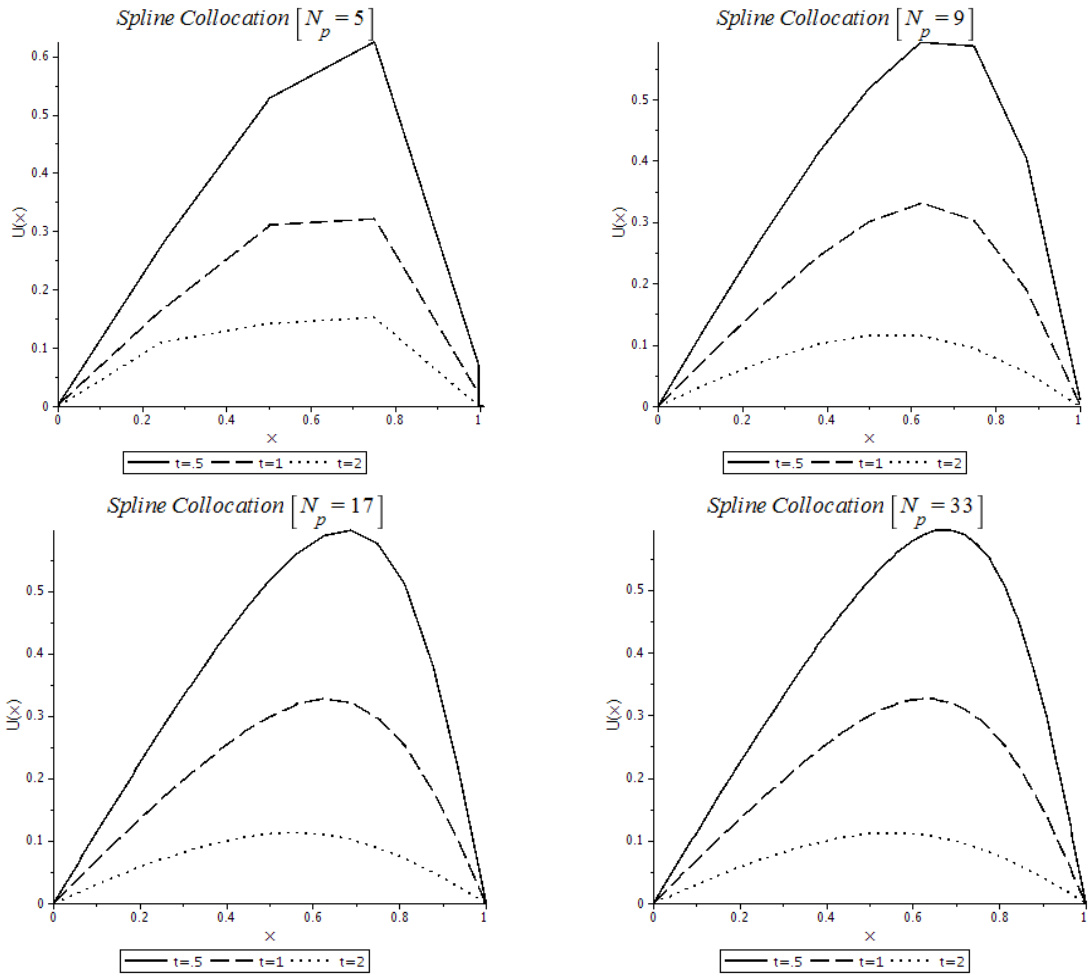

**Figure 5.** Numerical solution in different scales for test case 2 with $\text{Re} = 10.0$

**Table 4.** Comparison of the results for test case 2 with $\text{Re} = 10.0$

| Method | time | x=.1 | x=.3 | x=.5 | x=.7 | x=.9 |
|---|---|---|---|---|---|---|
| IFDM[6] | 0.5 | 0.11328 | 0.33168 | 0.51713 | 0.59382 | 0.32153 |
| BEM[7] |  | 0.11263 | 0.32982 | 0.51499 | 0.59230 | 0.32011 |
| EXACT |  | 0.11266 | 0.33010 | 0.51540 | 0.59304 | 0.32175 |
| Current work ($N_p$=33) |  | 0.11267 | 0.33000 | 0.51544 | 0.59262 | 0.32174 |
| IFDM[6] | 1 | 0.06810 | 0.19819 | 0.30100 | 0.31966 | 0.15268 |
| BEM[7] |  | 0.06766 | 0.19632 | 0.29782 | 0.31555 | 0.14993 |
| EXACT |  | 0.06750 | 0.19647 | 0.29834 | 0.31656 | 0.15097 |
| Current work ($N_p$=33) |  | 0.06750 | 0.19639 | 0.29837 | 0.31623 | 0.15096 |
| IFDM[6] | 2 | 0.02997 | 0.08301 | 0.11326 | 0.10227 | 0.04209 |
| BEM[7] |  | 0.02969 | 0.08109 | 0.11003 | 0.09894 | 0.04077 |
| EXACT |  | 0.02929 | 0.08101 | 0.11020 | 0.09915 | 0.04070 |
| Current work ($N_p$=33) |  | 0.02929 | 0.08096 | 0.11024 | 0.09906 | 0.04069 |

**Table 5.** Average relative error for test case 2 with $Re = 1$ and $Re = 10$

| Method | $Re = 1$ | | | $Re = 10$ | | |
|---|---|---|---|---|---|---|
| | $t = 0.05$ | $t = 0.1$ | $t = 0.2$ | $t = 0.5$ | $t = 1.0$ | $t = 2.0$ |
| Current work | $4.84 \times 10^{-4}$ | $4.41 \times 10^{-4}$ | $4.30 \times 10^{-4}$ | $2.42 \times 10^{-4}$ | $3.23 \times 10^{-4}$ | $4.27 \times 10^{-4}$ |
| IFDM [6] | $1.34 \times 10^{-3}$ | $4.98 \times 10^{-3}$ | $1.97 \times 10^{-2}$ | $3.13 \times 10^{-3}$ | $9.53 \times 10^{-3}$ | $2.83 \times 10^{-2}$ |
| BEM [7] | $1.69 \times 10^{-3}$ | $2.47 \times 10^{-3}$ | $3.62 \times 10^{-3}$ | $1.65 \times 10^{-3}$ | $2.99 \times 10^{-3}$ | $4.00 \times 10^{-3}$ |

As it is evident from Table 5, the average relative error for the current work is an order of magnitude lower than the IFDM and the BEM results, despite lower number of collocation points ($N_p = 5, 9, 17, 33$) used relative to the IFDM (100 nodes).

### 3.3 Test case 3.

In this test case the Burgers equation is solved with Neumann boundary condition. The Re number is set to 10.

I.C.  $u(x,0) = 50(\frac{1}{2} - x)^3$  (28)

B.C.  $\frac{\partial u}{\partial x}(0,t) = \frac{\partial u}{\partial x}(1,t) = 0$

This test case shows the performance of the spline wavelet collocation method discussed earlier when Neumann boundary condition is specified. The solution obtained is given plotted in Fig. 6 for $N_p = 17$ and $N_p = 33$ at three different time steps.

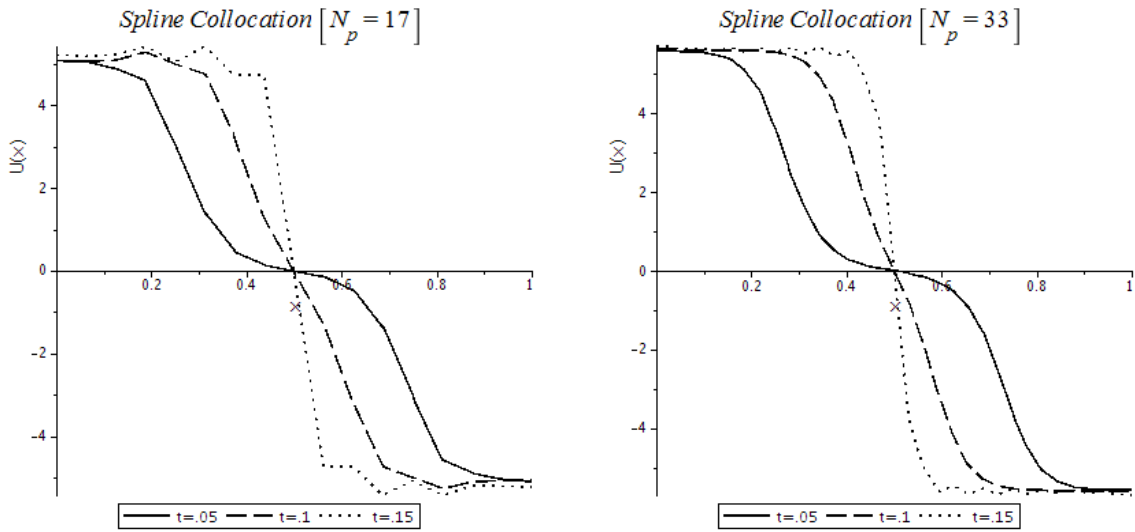

**Figure 6.** Numerical solutions in different scales for problem 3. And $Re = 10.0$

The solution of this problem is two fronts symmetrically advancing from the boundaries to the center of the domain as the time elapses [8]. This behavior is clearly seen in the numerical solution obtained in the current work, presented in Fig. 6. The solution is similar for $N_p = 17$ and $N_p = 33$ and the results are comparable; however, increasing the number of collocation points, makes the solution smoother. As the fronts advance and develops to a shock, wiggles appear around it (Fig 6-left). This is due to the fact that the spatial resolution of the solution is limited by number of basis functions and the smaller scales near the front are not resolved. By increasing the resolution of the numerical method, more scales are resolved, and the wiggles are diminished (Fig. 6-right). This shows one of the strengths of multi-resolution numerical schemes by which the details of the smaller scales can easily be added to the larger scale solution to obtain a refined solution.

## 4 Conclusion

In this paper a multi-resolution numerical method for solution of the PDEs is proposed and tested in different cases of Burgers equation. According to this method, the spatial operator is approximated by the compactly supported semi-orthogonal B-spline wavelets and the temporal operator is approximated by the second order finite difference method. Matrix of derivatives is used for expansion of derivatives terms. The obtained numerical results are compared with the existing exact and numerical solutions. Based on the test cases, the wavelet-based collocation method, gives better results than the fully implicit finite difference (IFDM) and mixed finite difference and boundary element method (BEM). The results obtained show that the relative error is an order of magnitude lower than the other methods, despite using lower number of collocation points.